\documentclass[12pt]{amsart}
\input{epsf}
\usepackage{amssymb,latexsym,cite}
\topskip  \numberwithin{equation}{section}
\newtheorem{theorem}{Theorem}[section]

\newtheorem{remark}[theorem]{Remark}

\begin{document}

\pagenumbering{arabic}
\pagestyle{headings}

\newcommand{\DPB}[4]{P\beta_{#1}^{(#2)}(#3,#4)}

\title[Ordered Bell and degenerate ordered Bell polynomials]{Representations by ordered Bell and degenerate ordered Bell polynomials}
\author{Dae San Kim}
\address{Department of Mathematics, Sogang University, Seoul 121-742, Republic of Korea}
\email{dskim@sogang.ac.kr}

\author{Taekyun Kim}
\address{Department of Mathematics, Kwangwoon University, Seoul 139-701, Republic of Korea}
\email{tkkim@kw.ac.kr}


\begin{abstract}
In this paper, we consider the problem of representing any polynomial in terms of the ordered Bell and degenerate ordered Bell polynomials, and more generally of the higher-order ordered Bell and higher-order degenerate ordered Bell polynomials. We derive explicit formulas with the help of umbral calculus and illustrate our results with some examples.
\end{abstract}
\maketitle


\section{Introduction and preliminaries}

The ordered Bell numbers (also called Fubini numbers) $b_n$ arise from number theory and various counting problems in enumerative combinatorics (see [5,14]). The ordered Bell numbers $b_n$ appeared already in 1859 work of Cayley [3], who used them to count certain plane trees with $n+1$ totally ordered leaves. While the (unordered) Bell numbers $Bel_{n}$ given by $e^{(e^t -1)}=\sum_{n=0}^{\infty}Bel_{n}\frac{t^n}{n!}$ count partitions of $[n]=\left\{1,2,\dots,n\right\}$ into nonempty disjoint subsets, the ordered Bell numbers count ordered partitions of $[n]$. Equivalently, the ordered Bell numbers $b_n$ count either the number of weak orderings on a set of $n$ elements or the mappings from $[n]$ to itself whose image is $[\ell]$, $1\leq \ell \le n$. They also count formulas in Fubini's theorem when rearranging the order of summation in multiple sums. We let the reader refer to [19], for details on the numerous uses of the ordered Bell numbers in counting problems.

Let $p(x) \in \mathbb{C}[x]$, with $\mathrm{deg}\, p(x)=n$. Write $p(x)=\sum_{k=0}^{n}a_kB_k(x)$, where $B_{n}(x)$ are the Bernoulli polynomials (see \eqref{3A}). Then it is known (see [12]) that
\begin{equation}\label{1A}
a_0=\int_{0}^{1}p(t) dt,\quad a_{k}=\frac{1}{k!}(p^{(k-1)}(1)-p^{(k-1)}(0)), \,\,\mathrm{for}\,\, k=1,2,\dots,n.
\end{equation}
Applying the formulas in \eqref{1A} to the polynomial $p(x)=\sum_{k=1}^{n-1}\frac{1}{k(n-k)}B_{k}(x)B_{n-k}(x)$, we can obtain an identity (see [12,16]) which yields, afer slight modification, the following identity:
\begin{align}
& \sum_{k=1}^{n-1}\frac{1}{2k\left(2n-2k\right)}B_{2k}\left(x\right)B_{2n-2k}\left(x\right)+\frac{2}{2n-1}B_{1}\left(x\right)B_{2n-1}\left(x\right)\label{2A}\\
= & \frac{1}{n}\sum_{k=1}^{n}\frac{1}{2k}\dbinom{2n}{2k}B_{2k}B_{2n-2k}\left(x\right)+\frac{1}{n}H_{2n-1}B_{2n}\left(x\right)\nonumber
+\frac{2}{2n-1}B_{1}\left(x\right)B_{2n-1},\nonumber
\end{align}
where $n \ge 2$, and $H_{n}=1+\frac{1}{2}+ \cdots +\frac{1}{n}$.
Letting $x=0$ and $x=\frac{1}{2}$ in \eqref{2A} respectively give a slight variant of the Miki's identity and the Faber-Pandharipande-Zagier (FPZ) identity.
Here it should be emphasized that the other proofs of Miki's (see [8,17,21]) and FPZ identities (see [6,7]) are quite involved, while our proofs of Miki's and FPZ identities follow from the simple formulas in \eqref{1A} involving only derivatives and integrals of the given polynomials. \par

Analogous formulas to Remark 3.2  can be obtained for the representations by Bernoulli, Euler and Genocchi polynomials. Many interesting identities have been derived by using these formulas (see [10,12,13]). The list in the References are far from being exhaustive. However, the interested reader can easily find more related papers in the literature. Also, we should mention here that there are other ways of obtaining the same result as the one in \eqref{2A}. One of them is to use Fourier series expansion of the function obtained by extending by periodicity 1 of the polynomial function restricted to the interval $[0,1)$ (see [15,16]).

The aim of this paper is to derive formulas (see Theorem 3.1) expressing any polynomial in terms of the degenerate ordered Bell polynomials (see [11]) with the help of umbral calculus (see [4,20,22]) and to illustrate our results with some examples (see Chapter 5). This is generalized to the higher-order degenerate ordered Bell polynomials. Indeed, we deduce formulas of representing any polynomial in terms of the higher-order degenerate ordered Bell polynomials (see Theorems 4.1) again by using umbral calculus and illustrate them with some examples (see Chapter 6). Letting $\lambda \rightarrow 0$, we obtain formulas of expressing any polynomial in terms of ordered Bell and higher-order ordered Bell polynomials. These formulas are also illustrated in Chapters 5 and 6. The contribution of this paper is the derivation of such formulas which have many potential applications. \par

The outline of this paper is as follows. In Section 1, we recall some necessary facts that are needed throughout this paper. In Section 2, we go over umbral calculus briefly. In Section 3, we derive formulas expressing any polynomial in terms of the degenerate ordered Bell polynomials. In Section 4, we derive formulas representing any polynomial in terms of the higher-order degenerate ordered Bell polynomials. In Section 5, we illustrate our results for representations by ordered Bell and higher-order ordered Bell polynomials with some examples. In Section 6, we illustrate our results for representations by degenerate ordered Bell and higher-order degenerate ordered Bell polynomials with some examples. Finally, we conclude our paper in Section 7.

The Bernoulli polynomials $B_n(x)$ are defined by
\begin{equation}\label{3A}
\frac{t}{e^t-1}e^{xt}=\sum_{n=0}^{\infty}B_{n}(x)\frac{t^n}{n!}.
\end{equation}
When $x=0$, $B_n=B_n(0)$ are called the Bernoulli numbers. We observe that $B_n(x)=\sum_{j=0}^{n}\binom{n}{j}B_{n-j}x^j,\,\frac{d}{dx}B_n(x)=nB_{n-1}(x),\,B_{n}(x+1)-B_{n}(x)=nx^{n-1}$.
The first few terms of $B_n$ are given by:
\begin{align*}
&B_0=1,\,B_1=-\frac{1}{2},\,B_2=\frac{1}{6},\,B_4=-\frac{1}{30},\,B_6=\frac{1}{42},\,B_8=-\frac{1}{30},\,B_{10}=\frac{5}{66},\\
&\,B_{12}=-\frac{691}{2730},\dots; B_{2k+1}=0,\,\,(k \ge 1).
\end{align*}

The Euler polynomials $E_n(x)$ are defined by
\begin{equation}\label{4A}
\frac{2}{e^t+1}e^{xt}=\sum_{n=0}^{\infty}E_{n}(x)\frac{t^n}{n!}.
\end{equation}
When $x=0$, $E_n=E_n(0)$ are called the Euler numbers. We observe that $E_n(x)=\sum_{j=0}^{n}\binom{n}{j}E_{n-j}x^j,\,\frac{d}{dx}E_n(x)=nE_{n-1}(x),\,E_{n}(x+1)+E_{n}(x)=2x^{n}$.
The first few terms of $E_n$ are given by:
\begin{align*}
&E_0=1,\,E_1=-\frac{1}{2},\,E_3=\frac{1}{4},\,E_5=-\frac{1}{2},\,E_7=\frac{17}{8},\,E_9=-\frac{31}{2},\dots;\\
&E_{2k}=0,\,\,(k \ge 1).
\end{align*}

The Genocchi polynomials $G_n(x)$ are defined by
\begin{equation}\label{5A}
\frac{2t}{e^t+1}e^{xt}=\sum_{n=0}^{\infty}G_{n}(x)\frac{t^n}{n!}.
\end{equation}
When $x=0$, $G_n=G_n(0)$ are called the Genocchi numbers. We observe that $G_n(x)=\sum_{j=0}^{n}\binom{n}{j}G_{n-j}x^j$,\,$\frac{d}{dx}G_n(x)=nG_{n-1}(x)$,\, $G_{n}(x+1)+G_{n}(x)=2nx^{n-1}$, and $\mathrm{deg}\,G_n(x)=n-1$, for $n \ge1$. Note that $G_{0}(x)=0$.
The first few terms of $G_n$ are given by:
\begin{align*}
&G_0=0,\,G_1=1,\,G_2=-1,\,G_4=1,\,G_6=-3,\,G_8=17,\,G_{10}=-155\\
&G_{12}=2073,\dots; G_{2k+1}=0,\,\,(k \ge 1).
\end{align*}

The ordered Bell polynomials $b_{n}(x)$ are defined by
\begin{equation}\label{6A}
\frac{1}{2-e^t}e^{xt}=\sum_{n=0}^{\infty}b_{n}(x)\frac{t^n}{n!}.
\end{equation}
When $x=0$, $b_n=b_n(0)$ are called the ordered Bell numbers. We observe that $b_n(x)=\sum_{j=0}^{n}\binom{n}{j}b_{n-j}x^j,\,\frac{d}{dx}b_n(x)=nb_{n-1}(x),\,2b_{n}(x)-b_{n}(x+1)=x^{n}$.
The first few terms of $b_n$ are given by:
\begin{align*}
b_0=1,\,b_1=1,\,b_2=3,\,b_3=13,\,b_4=75,\,b_5=541,\,b_6=4683,\,b_7=47293\dots.
\end{align*}

More generally, for any nonnegative integer $r$, the ordered Bell polynomials $b_n^{(r)}(x)$ of order $r$
are given by

\begin{equation}\label{7A}
\bigg(\frac{1}{2-e^t}\bigg)^{r}e^{xt}=\sum_{n=0}^{\infty}b_{n}^{(r)}(x)\frac{t^n}{n!}.
\end{equation}

For any nonzero real number $\lambda$, the degenerate exponentials are given by
\begin{align}
&e_{\lambda}^{x}(t)=(1+\lambda t)^{\frac{x}{\lambda}}=\sum_{n=0}^{\infty}(x)_{n,\lambda}\frac{t^n}{n!}, \label{8A}\\
&e_{\lambda}(t)=e_{\lambda}^{1}(t)=(1+\lambda t)^{\frac{1}{\lambda}}=\sum_{n=0}^{\infty}(1)_{n,\lambda}\frac{t^n}{n!}\nonumber.
\end{align}

In [11], in the spirit of [1] and as a degenerate version of them, the degenerate ordered Bell polynomials $b_{n,\lambda}(x)$ are introduced, which are given by
\begin{equation}\label{9A}
\frac{1}{2- e_{\lambda}(t)}e_{\lambda}^{x}(t)=\sum_{n=0}^{\infty}b_{n,\lambda}(x)\frac{t^n}{n!}.
\end{equation}
For $x=0$, $b_{n,\lambda}=b_{n,\lambda}(0)$ are called the degenerate ordered Bell numbers. \par
More generally, for any nonnegative integer $r$, the degenerate ordered Bell polynomials $b_{n,\lambda}^{(r)}(x)$ of order $r$
are given by

\begin{equation}\label{10A}
\bigg(\frac{1}{2-e_{\lambda}(t)}\bigg)^{r}e_{\lambda}^{x}(t)=\sum_{n=0}^{\infty}b_{n,\lambda}^{(r)}(x)\frac{t^n}{n!}.
\end{equation}

We remark that $b_{n,\lambda}(x) \rightarrow b_{n}(x)$, and $b_{n,\lambda}^{(r)}(x) \rightarrow b_{n}^{(r)}(x)$, as $\lambda$ tends to $0$.

We recall some notations and facts about forward differences. Let $f$ be any complex-valued function of the real variable $x$. Then, for any real number $a$, the forward difference $\Delta_{a}$ is given by
\begin{equation}\label{11A}
\Delta_{a}f(x)=f(x+a)-f(x).
\end{equation}
If $a=1$, then we let
\begin{equation}\label{12A}
\Delta f(x)=\Delta_{1}f(x)=f(x+1)-f(x).
\end{equation}

In general, the $n$th oder forward differences are given by
\begin{equation}\label{13A}
\Delta_{a}^{n}f(x)=\sum_{i=0}^{n}\binom{n}{i} (-1)^{n-i}f(x+ia).
\end{equation}
For $a=1$, we have
\begin{equation}\label{14A}
\Delta^{n}f(x)=\sum_{i=0}^{n}\binom{n}{i} (-1)^{n-i}f(x+i).
\end{equation}
Finally, we recall that the Stirling numbers of the second kind $S_{2}(n,k)$ can be given by means of
\begin{equation}\label{15A}
\frac{1}{k!}(e^{t}-1)^{k}=\sum_{n=k}^{\infty}S_{2}(n,k)\frac{t^{n}}{n!}.
\end{equation}

\section{Review of umbral calculus}
\vspace{0.5cm}
Here we will briefly go over very basic facts about umbral calculus. For more details on this, we recommend the reader to refer to [4,20,22].
Let $\mathbb{C}$ be the field of complex numbers. Then $\mathcal{F}$ denotes the algebra of formal power series in $t$ over $\mathbb{C}$, given by
\begin{displaymath}
 \mathcal{F}=\bigg\{f(t)=\sum_{k=0}^{\infty}a_{k}\frac{t^{k}}{k!}~\bigg|~a_{k}\in\mathbb{C}\bigg\},
\end{displaymath}
and $\mathbb{P}=\mathbb{C}[x]$ indicates the algebra of polynomials in $x$ with coefficients in $\mathbb{C}$. \par
Let $\mathbb{P}^{*}$ be the vector space of all linear functionals on $\mathbb{P}$. If $\langle L|p(x)\rangle$ denotes the action of the linear functional $L$ on the polynomial $p(x)$, then the vector space operations on $\mathbb{P}^{*}$ are defined by
\begin{displaymath}
\langle L+M|p(x)\rangle=\langle L|p(x)\rangle+\langle M|p(x)\rangle,\quad\langle cL|p(x)\rangle=c\langle L|p(x)\rangle,
\end{displaymath}
where $c$ is a complex number. \par
For $f(t)\in\mathcal{F}$ with $\displaystyle f(t)=\sum_{k=0}^{\infty}a_{k}\frac{t^{k}}{k!}\displaystyle$, we define the linear functional on $\mathbb{P}$ by
\begin{equation}\label{1B}
\langle f(t)|x^{k}\rangle=a_{k}.
\end{equation}
From \eqref{1B}, we note that
\begin{equation*}
 \langle t^{k}|x^{n}\rangle=n!\delta_{n,k},\quad(n,k\ge 0),
\end{equation*}
where $\delta_{n,k}$ is the Kronecker's symbol. \par
Some remarkable linear functionals are as follows:
\begin{align}
&\langle e^{yt}|p(x) \rangle=p(y), \nonumber \\
&\langle e^{yt}-1|p(x) \rangle=p(y)-p(0), \label{2B} \\
& \bigg\langle \frac{e^{yt}-1}{t}\bigg |p(x) \bigg\rangle = \int_{0}^{y}p(u) du.\nonumber
\end{align}
Let
\begin{equation}\label{3B}
 f_{L}(t)=\sum_{k=0}^{\infty}\langle L|x^{k}\rangle\frac{t^{k}}{k!}.
\end{equation}
Then, by \eqref{1B} and \eqref{3B}, we get
\begin{displaymath}
    \langle f_{L}(t)|x^{n}\rangle=\langle L|x^{n}\rangle.
\end{displaymath}
That is, $f_{L}(t)=L$. Additionally, the map $L\longmapsto f_{L}(t)$ is a vector space isomorphism from $\mathbb{P}^{*}$ onto $\mathcal{F}$.\par  Henceforth, $\mathcal{F}$ denotes both the algebra of formal power series  in $t$ and the vector space of all linear functionals on $\mathbb{P}$. $\mathcal{F}$ is called the umbral algebra and the umbral calculus is the study of umbral algebra.
For each nonnegative integer $k$, the differential operator $t^k$ on $\mathbb{P}$ is defined by
\begin{equation}\label{4B}
t^{k}x^n=\left\{\begin{array}{cc}
(n)_{k}x^{n-k}, & \textrm{if $k\le n$,}\\
0, & \textrm{if $k>n$.}
\end{array}\right.
\end{equation}
Extending \eqref{4B} linearly, any power series
\begin{displaymath}
 f(t)=\sum_{k=0}^{\infty}\frac{a_{k}}{k!}t^{k}\in\mathcal{F}
\end{displaymath}
gives the differential operator on $\mathbb{P}$ defined by
\begin{equation}\label{5B}
 f(t)x^n=\sum_{k=0}^{n}\binom{n}{k}a_{k}x^{n-k},\quad(n\ge 0).
\end{equation}
It should be observed that, for any formal power series $f(t)$ and any polynomial $p(x)$, we have
\begin{equation}\label{6B}
\langle f(t) | p(x) \rangle =\langle 1 | f(t)p(x) \rangle =f(t)p(x)|_{x=0}.
\end{equation}
Here we note that an element $f(t)$ of $\mathcal{F}$ is a formal power series, a linear functional and a differential  operator. Some notable differential operators are as follows:
\begin{align}
&e^{yt}p(x)=p(x+y), \nonumber\\
&(e^{yt}-1)p(x)=p(x+y)-p(x), \label{7B}\\
&\frac{e^{yt}-1}{t}p(x)=\int_{x}^{x+y}p(u) du.\nonumber
\end{align}

The order $o(f(t))$ of the power series $f(t)(\ne 0)$ is the smallest integer for which $a_{k}$ does not vanish. If $o(f(t))=0$, then $f(t)$ is called an invertible series. If $o(f(t))=1$, then $f(t)$ is called a delta series. \par
For $f(t),g(t)\in\mathcal{F}$ with $o(f(t))=1$ and $o(g(t))=0$, there exists a unique sequence $s_{n}(x)$ (deg\,$s_{n}(x)=n$) of polynomials such that
\begin{equation} \label{8B}
\big\langle g(t)f(t)^{k}|s_{n}(x)\big\rangle=n!\delta_{n,k},\quad(n,k\ge 0).
\end{equation}
The sequence $s_{n}(x)$ is said to be the Sheffer sequence for $(g(t),f(t))$, which is denoted by $s_{n}(x)\sim (g(t),f(t))$. We observe from \eqref{8B} that
\begin{equation}\label{9B}
s_{n}(x)=\frac{1}{g(t)}p_{n}(x),
\end{equation}
where $p_{n}(x)=g(t)s_{n}(x) \sim (1,f(t))$.\par
In particular, if $s_{n}(x) \sim (g(t),t)$, then $p_{n}(x)=x^n$, and hence
\begin{equation}\label{10B}
s_{n}(x)=\frac{1}{g(t)}x^n.
\end{equation}

It is well known that $s_{n}(x)\sim (g(t),f(t))$ if and only if
\begin{equation}\label{11B}
\frac{1}{g\big(\overline{f}(t)\big)}e^{x\overline{f}(t)}=\sum_{k=0}^{\infty}\frac{s_{k}(x)}{k!}t^{k},
\end{equation}
for all $x\in\mathbb{C}$, where $\overline{f}(t)$ is the compositional inverse of $f(t)$ such that $\overline{f}(f(t))=f(\overline{f}(t))=t$. \par

The following equations \eqref{12B}, \eqref{13B}, and \eqref{14B} are equivalent to the fact that  $s_{n}\left(x\right)$ is Sheffer for $\left(g\left(t\right),f\left(t\right)\right)$, for some invertible $g(t)$:
\begin{align}
f\left(t\right)s_{n}\left(x\right)&=ns_{n-1}\left(x\right),\quad\left(n\ge0\right),\label{12B}\\
s_{n}\left(x+y\right)&=\sum_{j=0}^{n}\binom{n}{j}s_{j}\left(x\right)p_{n-j}\left(y\right),\label{13B}
\end{align}
with $p_{n}\left(x\right)=g\left(t\right)s_{n}\left(x\right),$
\begin{equation}\label{14B}
s_{n}\left(x\right)=\sum_{j=0}^{n}\frac{1}{j!}\big\langle{g\left(\overline{f}\left(t\right)\right)^{-1}\overline{f}\left(t\right)^{j}}\big |{x^{n}\big\rangle}x^{j}.
\end{equation}

For $s_{n}(x)\sim(g(t),f(t))$, and $r_{n}(x)\sim(h(t),l(t))$, we have
\begin{equation}\label{15B}
s_{n}\left(x\right)=\sum_{k=0}^{n}C_{n,k}r_{k}\left(x\right),\quad\left(n\ge0\right),
\end{equation}
where
\begin{equation}\label{16B}
C_{n,k}=\frac{1}{k!}\bigg\langle{\frac{h\left(\overline{f}\left(t\right)\right)}{g\left(\overline{f}\left(t\right)\right)}l\left(\overline{f}\left(t\right)\right)^{k}}\Big |{x^{n}}\bigg \rangle.
\end{equation}

\section{Representation by degenerate ordered Bell polynomials}

Our interest here is to derive formulas expressing any polynomial in terms of the degenerate ordered Bell polynomials.

From \eqref{9A} and \eqref{11B}, we first observe that
\begin{align}
b_{n,\lambda}(x) \sim \big(g(t) &=2-e^t , f(t)=\frac{1}{\lambda}(e^{\lambda t}-1)\big),\label{1C} \\
& (x)_{n,\lambda} \sim (1, f(t)=\frac{1}{\lambda}(e^{\lambda t}-1)).\label{2C}
\end{align}
From \eqref{11A}, \eqref{7B}, \eqref{12B}, \eqref{1C} and \eqref{2C}, we note that
\begin{align}
f(t)b_{n,\lambda}(x)&=nb_{n-1,\lambda}(x)=\frac{1}{\lambda}(e^{\lambda t}-1)b_{n}(x)=\frac{1}{\lambda}\Delta_{\lambda}b_{n}(x), \label{3C}\\
&f(t)(x)_{n,\lambda}=n(x)_{n-1,\lambda}.\label{4C}
\end{align}
It is immediate to see from \eqref{9A} and \eqref{12A} that
\begin{align}
2 b_{n,\lambda}(x)-b_{n,\lambda}(x+1)=b_{n,\lambda}(x)-\Delta b_{n,\lambda}(x)=(I-\Delta)b_{n,\lambda}(x) = (x)_{n,\lambda},\label{5C}
\end{align}
where $I$ is the identity operator. \par
Now, we assume that $p(x) \in \mathbb{C}[x]$ has degree $n$, and write $p(x)=\sum_{k=0}^{n}a_{k}b_{k,\lambda}(x)$.
Let $h(x)=2p(x)-p(x+1)$. Then, from \eqref{5C}, we have
\begin{align}
h(x)&=\sum_{k=0}^{n}a_{k}(2 b_{k,\lambda}(x)-b_{k,\lambda}(x+1))\label{6C} \\
&=\sum_{k=0}^{n}a_{k}(x)_{k,\lambda}. \nonumber
\end{align}
For $k \ge 0$, from \eqref{4C} and \eqref{6C} we obtain
\begin{align}
(f(t))^{k}h(x)&=(f(t))^{k}\sum_{l=0}^{n}a_l(x)_{l,\lambda} \label{7C} \\
&=\sum_{l=k}^{n}l(l-1) \cdots (l-k+1)a_{l}(x)_{l-k,\lambda}.\nonumber
\end{align}
Letting $x=0$ in \eqref{7C}, we finally get
\begin{equation}\label{8C}
a_{k}=\frac{1}{k!}(f(t))^{k}h(x)|_{x=0}=\frac{1}{k!}\langle (f(t))^{k}|h(x) \rangle,\,\,(k \ge 0).
\end{equation}
An alternative expression of \eqref{8C} is given by\\
\begin{align}
a_k&=\frac{1}{k! \lambda^{k}}\Delta_{\lambda}^{k}h(x)|_{x=0} \nonumber\\
&=\frac{1}{k! \lambda^{k}}(2\Delta_{\lambda}^{k}p(x)-\Delta_{\lambda}^{k}p(x+1))|_{x=0}\label{9C} \\
&=\frac{1}{k! \lambda^{k}}(2\Delta_{\lambda}^{k }p(0)-\Delta_{\lambda}^{k}p(1)),\nonumber
\end{align}
as $f(t)h(x)=\frac{1}{\lambda}(e^{\lambda t}-1)h(x)=\frac{1}{\lambda} \Delta_{\lambda}h(x)$.

From \eqref{13A}, we have another alternative expression of \eqref{8C} which is given by
\begin{align}
a_k&=\frac{1}{k! \lambda^{k}}\Delta_{\lambda}^{k}h(x)|_{x=0} \nonumber\\
&=\frac{1}{k! \lambda^{k}}\sum_{j=0}^{k}\binom{k}{j}(-1)^{k-j}(h(x+j \lambda)|_{x=0} \label{10C }\\
&=\frac{1}{k! \lambda^{k}}\sum_{j=0}^{k}\binom{k}{j}(-1)^{k-j}\big(2p(x+j \lambda)-p(x+1+j \lambda)\big)|_{x=0} \nonumber \\
&=\frac{1}{k! \lambda^{k}}\sum_{j=0}^{k}\binom{k}{j}(-1)^{k-j}\big(2p(j \lambda)-p(1+k \lambda)\big).\nonumber
\end{align}

By using \eqref{15A}, we obtain yet another expression of \eqref{8C} that is given by
\begin{align}
a_{k}&=\frac{1}{k!}\langle (f(t))^{k}|h(x) \rangle \nonumber\\
&=\frac{1}{\lambda^{k}}\bigg\langle \frac{1}{k!} (e^{\lambda t}-1)^{k}\bigg|h(x) \bigg\rangle \label{11C}\\
&=\frac{1}{\lambda^{k}}\bigg\langle \sum_{l=k}^{\infty}S_{2}(l,k)\frac{\lambda^{l}t^{l}}{l!}\bigg|h(x) \bigg\rangle \nonumber \\
&=\sum_{l=k}^{n}S_{2}(l,k)\frac{\lambda^{l-k}}{l!}\big(2p^{(l)}(0)-p^{(l)}(1)\big), \nonumber
\end{align}
where $p^{(l)}(x)=(\frac{d}{dx})^{l}p(x)$.

Finally, from \eqref{8C}--\eqref{11C}, we get the following theorem.

\begin{theorem}
Let $p(x) \in \mathbb{C}[x], \mathrm{deg}\, p(x)=n$. Let $f(t)=\frac{1}{\lambda}(e^{\lambda t}-1)$. Then we have
$p(x)=\sum_{k=0}^{n}a_k b_{k,\lambda}(x)$, \\
where
\begin{align*}
a_{k}&=\frac{1}{k!}(f(t))^{k}(2p(x)-p(x+1))|_{x=0} \\
&=\frac{1}{k!}\langle (f(t))^{k}|2p(x)-p(x+1) \rangle \\
&=\frac{1}{k! \lambda^{k}}\big \langle \big(e^{\lambda t}-1\big)^{k} \big |2p(x)-p(x+1)\big \rangle \\
&=\frac{1}{k! \lambda^{k}}(2\Delta_{\lambda}^{k}p(x)-\Delta_{\lambda}^{k}p(x+1))|_{x=0} \\
&=\frac{1}{k! \lambda^{k}}\sum_{j=0}^{k}\binom{k}{j}(-1)^{k-j}\big(2p(j \lambda)-p(1+j \lambda)\big) \\
&=\sum_{l=k}^{n}S_{2}(l,k)\frac{\lambda^{l-k}}{l!}\big(2p^{(l)}(0)-p^{(l)}(1)\big), \nonumber
\,\,\mathrm{for}\,\, k=0,1,\dots,n.
\end{align*}
\end{theorem}

\begin{remark}
Let $p(x) \in \mathbb{C}[x], with \,\,\mathrm{deg}\, p(x)=n$. Write $p(x)=\sum_{k=0}^{n}a_k b_k(x)$. As $\lambda$ tends to $0$, $f(t) \rightarrow t$. Thus we obtain the following result:
\begin{align}
a_{k}=\frac{1}{k!}(2p^{(k)}(0)-p^{(k)}(1)), \,\,\mathrm{for}\,\, k=0,1,\dots,n. \label{12C}
\end{align}
\end{remark}

\begin{remark}
Formulas similar to \eqref{12C} for Bernoulli, Euler and Genocchi polynomials have been applied to many polynomials in oder to obtain interesting identities for certain special polynomials and numbers \textnormal{(see [10,12,13])}. Some of the polynomials that have been considered are as follows:

(a) \begin{equation*}
\sum B_{i_1}(x) \cdots B_{i_r}(x)E_{j_1}(x) \cdots E_{j_s}(x)G_{k_1+1}(x) \cdots G_{k_t+1}(x)x^{l},
\end{equation*}
where the sum is over all nonnegative integers $ i_1,\cdots, i_r, j_1,\cdots,j_s,k_1,\cdots,k_t,l$ such that
$i_1+ \cdots+ i_r +j_1+\cdots+j_s+k_1+ \cdots+k_t+l=n$, and $r, s, t, l$ are nonnegative integers with $r+s+t \ge 1$. \\
(b) \begin{equation*}
\sum \frac{B_{i_1}(x) \cdots B_{i_r}(x)E_{j_1}(x) \cdots E_{j_s}(x)G_{k_1+1}(x) \cdots G_{k_t+1}(x)x^{l}}{i_1! \cdots i_r! j_1! \cdots j_s! (k_1+1)! \cdots(k_t+1)! l!}, \\
\end{equation*}
where the sum is over all nonnegative integers $ i_1,\cdots, i_r, j_1,\cdots,j_s,k_1,\cdots,k_t,l $ such that
$i_1+ \cdots+ i_r +j_1+\cdots+j_s+k_1+ \cdots+k_t+l=n$, and $r, s, t, l$ are nonnegative integers with $r+s+t \ge 1$. \\
(c) \begin{equation*}
\sum \frac{B_{i_1}(x) \cdots B_{i_r}(x)E_{j_1}(x) \cdots E_{j_s}(x)G_{k_1+1}(x) \cdots G_{k_t+1}(x)x^{l}}{i_1 \cdots i_r j_1 \cdots j_s (k_1+1) \cdots(k_t+1) l}, \\
\end{equation*}
where the sum is over all positive integers $ i_1,\cdots, i_r, j_1,\cdots,j_s,l$ and nonnegative integers $k_1,\cdots,k_t$ such that $i_1+ \cdots+ i_r +j_1+\cdots+j_s+k_1+ \cdots+k_t+l=n$, and $r, s, t, l$ are nonnegative integers with $r+s+t \ge 1$.
\end{remark}

\section{Representation by higher-order degenerate ordered Bell polynomials}

Our interest here is to derive formulas expressing any polynomial in terms of the higher-order degenerate ordered Bell polynomials. \par
With $g(t)=2-e^{t}, f(t)=\frac{1}{\lambda}(e^{\lambda t}-1)$, and from \eqref{10A} and \eqref{11B}, we note that
\begin{equation}\label{1D}
b_{n,\lambda}^{(r)}(x) \sim (g(t)^{r}, f(t)).
\end{equation}
Also, from \eqref{1D} and \eqref{12B}, we have
\begin{equation}\label{2D}
f(t)b_{n,\lambda}^{(r)}(x)=nb_{n-1,\lambda}^{(r)}(x),
\end{equation}
and from \eqref{10A}, it is immediate to see that
\begin{equation}\label{3D}
g(t)b_{n,\lambda}^{(r)}(x)=2b_{n,\lambda}^{(r)}(x)-b_{n,\lambda}^{(r)}(x+1)=b_{n,\lambda}^{(r-1)}(x).
\end{equation}
Thus, from \eqref{3D} we have $g(t)^{r}b_{n,\lambda}^{(r)}(x)=b_{n,\lambda}^{(0)}(x)=(x)_{n,\lambda}$. \par
Now, we assume that $p(x) \in \mathbb{C}[x]$ has degree $n$, and write $p(x)=\sum_{k=0}^{n}a_{k}b_{k,\lambda}^{(r)}(x)$. Then we have
\begin{equation}\label{4D}
g(t)^{r}p(x)=\sum_{l=0}^{n}a_{l}\,g(t)^{r}b_{l,\lambda}^{(r)}(x)=\sum_{l=0}^{n}a_{l}(x)_{l,\lambda}.
\end{equation}
By using \eqref{4C} and \eqref{4D}, we observe that
\begin{align}
f(t)^{k}g(t)^{r}p(x)&=\sum_{l=0}^{n}a_{l}f(t)^{k}(x)_{l,\lambda}\label{5D}\\
&=\sum_{l=k}^{n}a_{l}\,l(l-1) \cdots (l-k+1)(x)_{l-k,\lambda}.\nonumber
\end{align}
By evaluating \eqref{5D} at $x=0$, we obtain
\begin{equation}\label{6D}
a_k=\frac{1}{k!}f(t)^{k}g(t)^{r}p(x)|_{x=0}=\frac{1}{k!}\langle f(t)^{k}g(t)^{r} | p(x) \rangle.
\end{equation}
This also follows from the observation $\langle g(t)^{r}f(t)^{k} | b_{l,\lambda}^{(r)}(x) \rangle=l!\,\delta_{l,k}.$ \par

We would like to find more explicit expressions for \eqref{6D}.

For this purpose, we first observe that
\begin{align}\label{7D}
g(t)^{r}p(x)&=(I-\Delta)^{r}p(x)=2^{r}\sum_{j=0}^{r}\binom{r}{j}\big(-\frac{1}{2}\big)^{j}p(x+j).
\end{align}

Several alternative expressions of \eqref{6D} follow from \eqref{7D} and \eqref{13A}, which are given by
\begin{align}
a_k&=\frac{1}{k!}g(t)^{r}f(t)^{k}p(x)|_{x=0}\nonumber\\
&=\frac{1}{k!\lambda^{k}}(I-{\Delta})^{r}\Delta_{\lambda}^{k}p(x)|_{x=0}\label{8D}\\
&=\frac{1}{k!\lambda^{k}}(I-\Delta)^{r-1}\Delta_{\lambda}^{k}(2p(x)-p(x+1))|_{x=0}\nonumber\\
&=\frac{2^{r}}{k!\lambda^{k}}\sum_{l=0}^{k}\sum_{j=0}^{r}(-1)^{k+j-l}\binom{k}{l}\binom{r}{j}\frac{1}{2^{j}}p(j+l \lambda)\nonumber
,\,\,(0 \le k \le n).
\end{align}

Next, from \eqref{15A} we observe that
\begin{align}\label{9D}
f(t)^{k}p(x)=\frac{k!}{\lambda^{k}}\frac{1}{k!}(e^{\lambda t}-1)^{k}p(x)
=\frac{k!}{\lambda^{k}}\sum_{l=k}^{n}S_{2}(l,k)\frac{\lambda^{l}}{l!}t^{l}p(x),
\end{align}

Combining \eqref{7D} and \eqref{9D}, we see that
\begin{align}
a_k&=\frac{1}{k!}f(t)^{k}g(t)^{r}p(x)|_{x=0}\nonumber\\
&=\frac{1}{\lambda^{k}}\sum_{l=k}^{n}S_{2}(l,k)\frac{\lambda^{l}}{l!}t^{l}g(t)^{r}p(x)|_{x=0}\label{10D}\\
&=\frac{1}{\lambda^{k}}\sum_{l=k}^{n}S_{2}(l,k)\frac{\lambda^{l}}{l!}t^{l}
2^{r}\sum_{j=0}^{r}\binom{r}{j}\big(-\frac{1}{2}\big)^{j}p(x+j)|_{x=0}\nonumber\\
&=\frac{2^{r}}{\lambda^{k}}\sum_{l=k}^{n}\sum_{j=0}^{r}
\binom{r}{j}(-\frac{1}{2})^{j}\frac{\lambda^{l}}{l!}S_{2}(l,k)p^{(l)}(j).\nonumber
\end{align}

Now, from \eqref{8D} and \eqref{10D}, we finally arrive at the following theorem.
\begin{theorem}
Let $p(x) \in \mathbb{C}[x], \mathrm{deg}\, p(x)=n$. Let  $g(t)=2-e^{t}, f(t)=\frac{1}{\lambda}(e^{\lambda t}-1)$. Then we have $p(x)=\sum_{k=0}^{n}a_k b_{k,\lambda}^{(r)}(x)$, where
\begin{align*}
a_{k} &=\frac{1}{k!}f(t)^{k}g(t)^{r}p(x)|_{x=0}\\
&=\frac{1}{k!}\langle f(t)^{k}g(t)^{r} | p(x) \rangle \\
&=\frac{1}{k!\lambda^{k}}(I-{\Delta})^{r}\Delta_{\lambda}^{k}p(x)|_{x=0}\\
&=\frac{1}{k!\lambda^{k}}(I-\Delta)^{r-1}\Delta_{\lambda}^{k}(2p(x)-p(x+1))|_{x=0}\\
&=\frac{2^{r}}{k!\lambda^{k}}\sum_{l=0}^{k}\sum_{j=0}^{r}(-1)^{k+j-l}\binom{k}{l}\binom{r}{j}\frac{1}{2^{j}}p(j+l \lambda)\\
&=2^{r}\sum_{l=k}^{n}\sum_{j=0}^{r}
\binom{r}{j}(-\frac{1}{2})^{j}\frac{\lambda^{l-k}}{l!}S_{2}(l,k)p^{(l)}(j),\,\,(0 \le k \le n).
\end{align*}
\end{theorem}

\begin{remark}
Let $p(x) \in \mathbb{C}[x], with \,\,\mathrm{deg}\, p(x)=n$. Write $p(x)=\sum_{k=0}^{n}a_kb_{k}^{(r)}(x)$. As $\lambda$ tends to $0$, $f(t) \rightarrow t $. Thus, from Theorem 4.1, we have the following result:\par
\begin{align}
a_{k}&=\frac{1}{k!}f(t)^{k}g(t)^{r}p(x)|_{x=0}\nonumber\\
&=\frac{2^{r}}{k!}\sum_{j=0}^{r}\binom{r}{j}(-\frac{1}{2})^{j}p^{(k)}(j)\nonumber\\
&=\frac{1}{k!}(I-\Delta)^{r}p^{(k)}(x)|_{x=0}\nonumber\\
&=\frac{1}{k!}\sum_{j=0}^{r}\binom{r}{j}(-1)^{j}\Delta^{j}p^{(k)}(x)|_{x=0}\label{ } \\
&=\frac{1}{k!}(I-\Delta)^{r-1}\big(2p^{(k)}(x)-p^{(k)}(x+1)\big)|_{x=0} \nonumber\\
&=\frac{1}{k!}\sum_{j=0}^{r-1}\binom{r-1}{j}(-1)^{j}\big(2\Delta^{j}p^{(k)}(x)- \Delta^{j}p^{(k)}(x+1)\big)|_{x=0},\nonumber
\end{align}
where $g(t)=2-e^{t},\,\,f(t)=t$.
\end{remark}

\section{Examples on representations by ordered Bell polynomials}

Here we illustrate our formulas in Remarks 3.2 and 4.2 with some examples.

(a) Let $p(x)=B_{n}(x)=\sum_{k=0}^{n}a_{k}b_{k}(x)$. Then, noting that $B_{n}(1)-B_{n}=\delta_{n,1}$ from \eqref{3A}, we have \par
\begin{align*}
a_{k}&=\frac{1}{k!}(2p^{(k)}(0)-p^{(k)}(1))=\binom{n}{k}(2B_{n-k}-B_{n-k}(1))\\
&=\binom{n}{k}\big(B_{n-k}-(B_{n-k}(1)-B_{n-k})\big) \\
&=\binom{n}{k}(B_{n-k}-\delta_{n-k,1}).
\end{align*}
Thus $B_{n}(x)=\sum_{k=0}^{n}\binom{n}{k}(B_{n-k}-\delta_{n-k,1})b_{k}(x)$. \par
Now, we let $B_{n}(x)=\sum_{k=0}^{n}c_{k}b_{k}^{(r)}(x)$. Then we obtain the following:
\begin{align*}
c_{k}&=\frac{2^{r}}{k!}\sum_{j=0}^{r}\binom{r}{j}\big(-\frac{1}{2}\big)^{j}p^{(k)}(j)\\
&=2^{r}\binom{n}{k}\sum_{j=0}^{r}\binom{r}{j}(-\frac{1}{2})^{j}B_{n-k}(j).
\end{align*}
Hence $B_{n}(x)=2^{r}\sum_{k=0}^{n}\sum_{j=0}^{r}\binom{n}{k}\binom{r}{j}\big(-\frac{1}{2}\big)^{j}B_{n-k}(j)\,b_{k}^{(r)}(x)$.

These results on representations of $B_{n}(x)$ by $b_{k}(x)$ and $b_{k}^{(r)}(x)$ also follow from \eqref{16B}.

(b) Let $p(x)=\sum_{k=1}^{n-1}\frac{1}{k(n-k)}B_{k}(x)B_{n-k}(x),\,\,(n \ge 2)$. For this, we first recall from [12] that
\begin{equation}\label{1E}
p(x)=\frac{2}{n}\sum_{l=0}^{n-2}\frac{1}{n-l}\binom{n}{l}B_{n-l}B_l(x)+\frac{2}{n}H_{n-1}B_n(x),
\end{equation}
where $H_{n}=1+\frac{1}{2}+\cdots+\frac{1}{n}$ is the harmonic number. A slight modification of \eqref{1E} gives the Miki's and FPZ identities.
Let $p(x)=\sum_{k=0}^{n}a_{k}b_{k}(x)$. Then, for  $0 \le k \le n-2$, we have
\begin{align*}
k!a_{k}&=2p^{(k)}(0)-p^{(k)}(1)\\
&=\frac{2}{n}\sum_{l=k}^{n-2}\frac{1}{n-l}\binom{n}{l}B_{n-l}(l)_{k}(B_{l-k}-\delta_{l-k,1})+\frac{2}{n}H_{n-1}(n)_{k}B_{n-k};
\end{align*}
For $k=n-1$ or $k=n$, we have \par
\begin{align*}
k!a_{k}=\frac{2}{n}H_{n-1}(n)_{k}(B_{n-k}-\delta_{n-k,1}).
\end{align*}
Hence we get \par
\begin{align*}
&\sum_{k=1}^{n-1}\frac{1}{k(n-k)}B_{k}(x)B_{n-k}(x)\\
&=\sum_{k=0}^{n-2}\bigg\{\frac{2}{n}\sum_{l=k}^{n-2}\frac{1}{n-l}\binom{n}{l}\binom{l}{k}B_{n-l}(B_{l-k}-\delta_{l-k,1})+\frac{2}{n}\binom{n}{k}H_{n-1}B_{n-k}\bigg\}b_{k}(x) \\
&\quad\quad\quad\quad-3H_{n-1}b_{n-1}(x)+\frac{2}{n}H_{n-1}b_{n}(x).
\end{align*}

(c) Let $\Omega_{l}=\sum_{a=1}^{s}\binom{s}{a}2^{a}(-1)^{s-a}\sum_{i_{1}+\cdots+i_{a}=l}\prod_{j=1}^{a}b_{i_{j}}-\sum_{i_{1}+\cdots+i_{s}=l}\prod_{j=1}^{s}b_{i_{j}}$.
Let $p(x)=\sum_{i_{1}+\cdots+i_{s}=n}\prod_{j=1}^{s}b_{i_{j}}(x)$. Then, from Theorems 2 and 3 in [5],  it is immediate to see that
\begin{align*}
p(x)=\frac{1}{n+s}\sum_{j=0}^{n}\binom{n+s}{j}\Omega_{n-j+1}B_{j}(x).
\end{align*}
Let $p(x)=\sum_{k=0}^{n}a_{k}b_{k}(x)$. Then, for $0 \le k \le n$, we have \par
\begin{align*}
k!a_{k}=2p^{(k)}(0)-p^{(k)}(1)=\frac{1}{n+s}\sum_{j=k}^{n}\binom{n+s}{j}\Omega_{n-j+1}(j)_{k}(B_{j-k}-\delta_{j-k,1}). \\
\end{align*}
Thus we obtain \par
\begin{align*}
\sum_{i_{1}+\cdots+i_{s}=n}\prod_{j=1}^{s}b_{i_{j}}(x)=\frac{1}{n+s}\sum_{k=0}^{n}\bigg\{\sum_{j=k}^{n}\binom{n+s}{j}\binom{j}{k}\Omega_{n-j+1}(B_{j-k}-\delta_{j-k,1})\bigg\}b_{k}(x).
\end{align*}

(d) In [16], it is proved that the following identity is valid for $n \ge 2$:
\begin{equation}\label{7E}
\sum_{k=1}^{n-1}\frac{1}{k(n-k)}G_k(x)G_{n-k}(x)=-\frac{4}{n}\sum_{k=0}^{n-2}\binom{n}{k}\frac{G_{n-k}}{n-k}B_k(x).
\end{equation}
Write $p(x)=\sum_{k=1}^{n-1}\frac{1}{k(n-k)}G_k(x)G_{n-k}(x)=\sum_{k=0}^{n-2}a_{k}b_{k}^{(r)}(x).$
\begin{align*}
a_{k}&=\frac{2^{r}}{k!}\sum_{j=0}^{r}\binom{r}{j}\big(-\frac{1}{2}\big)^{j}p^{(k)}(j)\\
&=-\frac{2^{r+2}}{n}\sum_{j=0}^{r}\sum_{l=k}^{n-2}\big(-\frac{1}{2}\big)^{j}\binom{r}{j}\binom{n}{l}\binom{l}{k}\frac{G_{n-l}}{n-l}B_{l-k}(j).
\end{align*}
Therefore we obtain the following result:
\begin{align*}
&\sum_{k=1}^{n-1}\frac{1}{k(n-k)}G_k(x)G_{n-k}(x)\\
&=-\frac{2^{r+2}}{n}\sum_{k=0}^{n-2}\bigg\{\sum_{j=0}^{r}\sum_{l=k}^{n-2}(-\frac{1}{2})^{j}\binom{r}{j}\binom{n}{l}\binom{l}{k}\frac{G_{n-l}}{n-l}B_{l-k}(j)\bigg\}b_{k}^{(r)}(x).
\end{align*}

(e) Nielsen [18,2] represented products of two Euler polynomials in terms of Bernoulli polynomials as follows:
\begin{align*}
E_{m}(x)E_{n}(x)&=-2\sum_{i=1}^{m}\binom{m}{i}E_{i} \frac{B_{m+n-i+1}(x)}{m+n-i+1}\\
&\quad-2\sum_{j=1}^{n}\binom{n}{j}E_{j} \frac{B_{m+n-j+1}(x)}{m+n-j+1}\\
&\quad +2(-1)^{n+1} \frac{m!\,n!}{(m+n+1)!} E_{m+n+1}.
\end{align*}
Write $p(x)=E_{m}(x)E_{n}(x)=\sum_{k=0}^{m+n}a_{k}b_{k}^{(r)}(x)$.
We now observe that
\begin{align*}
a_{k}&=\frac{2^{r}}{k!}\sum_{j=0}^{r}\binom{r}{j}\big(-\frac{1}{2}\big)^{j}p^{(k)}(j)\\
&=\frac{2^{r}}{k!}\sum_{j=0}^{r}\sum_{a+b=k}\binom{r}{j}\binom{k}{a,b}\big(-\frac{1}{2}\big)^{j}(E_{m}(x))^{(a)}(E_{n}(x))^{(b)}|_{x=j},
\end{align*}
where $(E_{n}(x))^{(a)}=(\frac{d}{dx})^{a}E_{n}(x)$, and $\binom{k}{a,b}=\frac{k!}{a!b!}$.
Thus we obtain the following:
\begin{align*}
-2\sum_{i=1}^{m}&\binom{m}{i}E_{i} \frac{B_{m+n-i+1}(x)}{m+n-i+1}
-2\sum_{j=1}^{n}\binom{n}{j}E_{j} \frac{B_{m+n-j+1}(x)}{m+n-j+1}\\
&+2(-1)^{n+1} \frac{m!\,n!}{(m+n+1)!} E_{m+n+1}\\
&=\sum_{k=0}^{m+n}\frac{2^{r}}{k!}\bigg\{\sum_{j=0}^{r}\sum_{a+b=k}\binom{r}{j}\binom{k}{a,b}(-\frac{1}{2})^{j}(E_{m}(x))^{(a)}(E_{n}(x))^{(b)}|_{x=j}\bigg\}b_{k}^{(r)}(x).
\end{align*}

\section{Examples on representations by degenerate ordered Bell polynomials}

Here we illustrate our formulas in Theorems 3.1 and 4.1 with some examples.

(a) Let $p(x)=E_{n}(x)=\sum_{k=0}^{n}a_{k}b_{k,\lambda}(x)$. Then, noting $E_{n}(1)+E_{n}=2\delta_{n,0}$, we have \par
\begin{align*}
a_{k}&=\frac{1}{k! \lambda^{k}}(2\Delta_{\lambda}^{k}E_{n}(x)-\Delta_{\lambda}^{k}E_{n}(x+1))|_{x=0} \\
&=\frac{1}{k! \lambda^{k}}\sum_{j=0}^{k}\binom{k}{j}(-1)^{k-j}\big(2E_{n}(j \lambda)-E_{n}(1+j \lambda)\big) \\
&=\sum_{l=k}^{n}\binom{n}{l}\lambda^{l-k}S_{2}(l,k)\big(3E_{n-l}-2\delta_{n,l}\big).
\end{align*}
Hence we obtain
\begin{align*}
E_{n}(x)&=\sum_{k=0}^{n}\frac{1}{k! \lambda^{k}}\big((2\Delta_{\lambda}^{k}E_{n}(x)-\Delta_{\lambda}^{k}E_{n}(x+1))|_{x=0}\big)b_{k,\lambda}(x)  \\
&=\sum_{k=0}^{n}\frac{1}{k! \lambda^{k}}\bigg\{\sum_{j=0}^{k}\binom{k}{j}(-1)^{k-j}\big(2E_{n}(j \lambda)-E_{n}(1+j \lambda)\big)\bigg\}b_{k,\lambda}(x) \\
&=\sum_{k=0}^{n}\bigg\{\sum_{l=k}^{n}\binom{n}{l}\lambda^{l-k}S_{2}(l,k)\big(3E_{n-l}-2\delta_{n,l}\big)\bigg\}b_{k,\lambda}(x).
\end{align*}

Now, we let $p(x)=E_{n}(x)=\sum_{k=0}^{n}c_{k}b_{k,\lambda}^{(r)}(x)$. Recalling that $E_{n}(x+1)+E_{n}(x)=2x^{n}$, we get
\begin{align*}
c_{k}&=\frac{1}{k!\lambda^{k}}(I-{\Delta})^{r}\Delta_{\lambda}^{k}E_{n}(x)|_{x=0}\\
&=\frac{1}{k!\lambda^{k}}(I-\Delta)^{r-1}\Delta_{\lambda}^{k}(3E_{n}(x)-2x^{n})|_{x=0}
\end{align*}
\begin{align*}
&=\frac{2^{r}}{k!\lambda^{k}}\sum_{l=0}^{k}\sum_{j=0}^{r}(-1)^{k+j-l}\binom{k}{l}\binom{r}{j}\frac{1}{2^{j}}E_{n}(j+l \lambda)\\
&=\frac{2^{r}}{\lambda^{k}}\sum_{l=k}^{n}\sum_{j=0}^{r}
\binom{r}{j}\binom{n}{l}(-\frac{1}{2})^{j}\lambda^{l}S_{2}(l,k)E_{n-l}(j).
\end{align*}

This shows the following:
\begin{align*}
E_{n}(x)&=\sum_{k=0}^{n}\frac{1}{k!\lambda^{k}}\big((I-{\Delta})^{r}\Delta_{\lambda}^{k}E_{n}(x)|_{x=0}\big)b_{k,\lambda}^{(r)}(x)\\
&=\sum_{k=0}^{n}\frac{1}{k!\lambda^{k}}\big((I-\Delta)^{r-1}\Delta_{\lambda}^{k}(3E_{n}(x)-2x^{n})|_{x=0}\big)b_{k,\lambda}^{(r)}(x)\\
&=\sum_{k=0}^{n}\bigg\{\frac{2^{r}}{k!\lambda^{k}}\sum_{l=0}^{k}\sum_{j=0}^{r}(-1)^{k+j-l}\binom{k}{l}\binom{r}{j}\frac{1}{2^{j}}E_{n}(j+l \lambda)\bigg\}b_{k,\lambda}^{(r)}(x) \\
&=\sum_{k=0}^{n}\bigg\{\frac{2^{r}}{\lambda^{k}}\sum_{l=k}^{n}\sum_{j=0}^{r}
\binom{r}{j}\binom{n}{l}(-\frac{1}{2})^{j}\lambda^{l}S_{2}(l,k)E_{n-l}(j)\bigg\}b_{k,\lambda}^{(r)}(x).
\end{align*}
These results on representations of $E_{n}(x)$ by $b_{k,\lambda}(x)$ and $b_{k,\lambda}^{(r)}(x)$ also follow from \eqref{16B}.

(b) Working similarly to (a) and recalling that $2b_{n}(x)-b_{n}(x+1)=x^{n}$, we have
\begin{align*}
b_{n}(x)&=\sum_{k=0}^{n}\bigg\{\frac{1}{k! \lambda^{k}}(2\Delta_{\lambda}^{k}b_{n}(x)-\Delta_{\lambda}^{k}b_{n}(x+1))|_{x=0}\bigg\}\,b_{k,\lambda}(x) \\
&=\sum_{k=0}^{n}\bigg\{\frac{1}{k! \lambda^{k}}\sum_{j=0}^{k}\binom{k}{j}(-1)^{k-j}\big(2b_{n}(j \lambda)-b_{n}(1+j \lambda)\big)\bigg\}\,b_{k,\lambda}(x)  \\
&=\sum_{k=0}^{n}S_{2}(n,k)\lambda^{n-k}b_{k,\lambda}(x).
\end{align*}

More generally, we also have

\begin{align*}
b_{n}(x)&=\sum_{k=0}^{n}\bigg\{\frac{1}{k!\lambda^{k}}(I-{\Delta})^{r}\Delta_{\lambda}^{k}b_{n}(x)|_{x=0}\bigg\}\,b_{k,\lambda}^{(r)}(x)\\
&=\sum_{k=0}^{n}\bigg\{\frac{1}{k!\lambda^{k}}(I-\Delta)^{r-1}\Delta_{\lambda}^{k}(2b_{n}(x)-b_{n}(x+1))|_{x=0}\bigg\}\,b_{k,\lambda}^{(r)}(x)\\
&=\sum_{k=0}^{n}\bigg\{\frac{2^{r}}{k!\lambda^{k}}\sum_{l=0}^{k}\sum_{j=0}^{r}(-1)^{k+j-l}\binom{k}{l}\binom{r}{j}\frac{1}{2^{j}}b_{n}(j+l \lambda)\bigg\}\,b_{k,\lambda}^{(r)}(x)\\
&=\sum_{k=0}^{n}\bigg\{\frac{2^{r}}{\lambda^{k}}\sum_{l=k}^{n}\sum_{j=0}^{r}
\binom{r}{j}\binom{n}{l}(-\frac{1}{2})^{j}\lambda^{l}S_{2}(l,k)b_{n-l}(j)\bigg\}\,b_{k,\lambda}^{(r)}(x).
\end{align*}

(c) As we mentioned earlier in \eqref{ }, the following identity holds:

\begin{align*}
\sum_{i_{1}+\cdots+i_{s}=n}\prod_{j=1}^{s}b_{i_{j}}(x)=\frac{1}{n+s}\sum_{i=0}^{n}\binom{n+s}{i}\Omega_{n-i+1}B_{i}(x),
\end{align*}
where $\Omega_{l}=\sum_{a=1}^{s}\binom{s}{a}2^{a}(-1)^{s-a}\sum_{i_{1}+\cdots+i_{a}=l}\prod_{j=1}^{a}b_{i_{j}}-\sum_{i_{1}+\cdots+i_{s}=l}\prod_{j=1}^{s}b_{i_{j}}$.
\par
Then, by proceeding similarly to (a), we have
\begin{align*}
&\sum_{i_{1}+\cdots+i_{s}=n}\prod_{j=1}^{s}b_{i_{j}}(x)\\
&=\sum_{k=0}^{n}\bigg\{\frac{1}{k!\lambda^{k}}\frac{1}{n+s}\sum_{i=0}^{n}\binom{n+s}{i}\Omega_{n-i+1}(I-{\Delta})^{r}\Delta_{\lambda}^{k}B_{i}(x)|_{x=0}\bigg\}b_{k,\lambda}^{(r)}(x)\\
&=\sum_{k=0}^{n}\bigg\{\frac{1}{n+s}\frac{1}{k!\lambda^{k}}\sum_{i=0}^{n}\binom{n+s}{i}\Omega_{n-i+1}(I-\Delta)^{r-1}\Delta_{\lambda}^{k}(B_{i}(x)-jx^{j-1})|_{x=0}\bigg\}b_{k,\lambda}^{(r)}(x)\\
&=\sum_{k=0}^{n}\bigg\{\frac{2^{r}}{k!\lambda^{k}}\frac{1}{n+s}\sum_{l=0}^{k}\sum_{j=0}^{r}\sum_{i=0}^{n}(-1)^{k+j-l}\binom{k}{l}\binom{r}{j}\binom{n+s}{i}\frac{1}{2^{j}}\Omega_{n-i+1}B_{i}(j+l \lambda)\bigg\}b_{k,\lambda}^{(r)}(x)\\
&=\sum_{k=0}^{n}\bigg\{\frac{2^{r}}{\lambda^{k}}\frac{1}{n+s}\sum_{l=k}^{n}\sum_{j=0}^{r}\sum_{i=l}^{n}
\binom{r}{j}\binom{i}{l}\binom{n+s}{i}(-\frac{1}{2})^{j}\lambda^{l}S_{2}(l,k)\Omega_{n-i+1}B_{i-l}(j)\bigg\}b_{k,\lambda}^{(r)}(x)
\end{align*}

(d) Nielsen [18,2] expressed products of two Bernoulli polynomials in terms of Bernoulli polynomials. Namely, for positive integers $m$ and $n$, with $m+n \ge 2$,
\begin{align*}
B_{m}(x)B_{n}(x)=\sum_{r}\left\{\binom{m}{2r}n+\binom{n}{2r}m\right\}\frac{B_{2r}B_{m+n-2r}(x)}{m+n-2r}+(-1)^{m+1}\frac{B_{m+n}}{\binom{m+n}{m}}.
\end{align*}
Then, by proceeding analogously to (a), we get
\begin{align*}
&\sum_{r}\left\{\binom{m}{2r}n+\binom{n}{2r}m\right\}\frac{B_{2r}B_{m+n-2r}(x)}{m+n-2r}+(-1)^{m+1}\frac{B_{m+n}}{\binom{m+n}{m}} \\
&=\sum_{k=0}^{m+n}\bigg\{\frac{1}{k!\lambda^{k}}(I-{\Delta})^{r}\Delta_{\lambda}^{k}\big(B_{m}(x)B_{n}(x)\big)|_{x=0}\bigg\}\, b_{k,\lambda}^{(r)}(x)\\
&=\sum_{k=0}^{m+n}\bigg\{\frac{1}{k!\lambda^{k}}(I-\Delta)^{r-1}\Delta_{\lambda}^{k}\big(2B_{m}(x)B_{n}(x)-B_{m}(x+1)B_{n}(x+1)\big)|_{x=0}\bigg\}\, b_{k,\lambda}^{(r)}(x)\\
&=\sum_{k=0}^{m+n}\bigg\{\frac{2^{r}}{k!\lambda^{k}}\sum_{l=0}^{k}\sum_{j=0}^{r}(-1)^{k+j-l}\binom{k}{l}\binom{r}{j}\frac{1}{2^{j}}B_{m}(j+l \lambda)B_{n}(j+l \lambda)\bigg\} \,b_{k,\lambda}^{(r)}(x)\\
&=\sum_{k=0}^{m+n}\bigg\{\frac{2^{r}}{\lambda^{k}}\sum_{l=k}^{n}\sum_{j=0}^{r}\sum_{a+b=l}
\binom{r}{j}\binom{l}{a,b}(-\frac{1}{2})^{j}\frac{\lambda^{l}}{l!}S_{2}(l,k)\big(B_{m}(x)\big)^{(a)}\big(B_{n}(x)\big)^{(b)}|_{x=j}\bigg\}\, b_{k,\lambda}^{(r)}(x),
\end{align*}
where $(B_{n}(x))^{(a)}=(\frac{d}{dx})^{a}B_{n}(x)$.

(e) In (4.20) of [9], it is shown that the following identity holds for $n \ge s$:
\begin{align*}
&\sum_{\substack{i_1+\cdots+i_s=n\\i_1,\dots,i_s \ge 1}} G_{i_{1}}(x) \dots G_{i_{s}}(x) \\
&\quad\quad=\frac{1}{n+s}\sum_{l=1}^{s}\binom{s}{l}(-2)^{l-1}\sum_{\substack{i_0 + i_1 + \cdots + i_{s-l} = n+1-l \\ i_0, i_1,\dots, i_{s-l} \ge 1}}\binom{n+s}{i_0} G_{i_1} \cdots G_{i_s}G_{i_0}(x).
\end{align*}

Then, proceeding similarly to (b), we can show that

\begin{align*}
&\sum_{\substack{i_1+\cdots+i_s=n\\i_1,\dots,i_s \ge 1}} G_{i_{1}}(x) \dots G_{i_{s}}(x)
=\sum_{k=0}^{n-s}\bigg\{\frac{1}{k!\lambda^{k}}\frac{1}{n+s}\sum_{l=1}^{s}\binom{s}{l}(-2)^{l-1}\\
&\quad \times\sum_{\substack{i_0 + i_1 + \cdots + i_{s-l} = n+1-l \\ i_0, i_1,\dots, i_{s-l} \ge 1}}\binom{n+s}{i_0} G_{i_1} \cdots G_{i_s}(I-{\Delta})^{r}\Delta_{\lambda}^{k}G_{i_0}(x)|_{x=0}\bigg\}\,b_{k,\lambda}^{(r)}(x)\\
&=\sum_{k=0}^{n-s}\bigg\{\frac{1}{k!\lambda^{k}}\frac{1}{n+s}\sum_{l=1}^{s}\binom{s}{l}(-2)^{l-1}\sum_{\substack{i_0 + i_1 + \cdots + i_{s-l} = n+1-l \\ i_0, i_1,\dots, i_{s-l} \ge 1}}\binom{n+s}{i_0} G_{i_1} \cdots G_{i_s}\\
&\quad \times (I-\Delta)^{r-1}\Delta_{\lambda}^{k}\big(3G_{i_0}(x)-2i_{0}x^{i_{0}-1}\big)|_{x=0}\bigg\}\,b_{k,\lambda}^{(r)}(x)\\
&=\sum_{k=0}^{n-s}\bigg\{\frac{2^{r}}{k!\lambda^{k}}\frac{1}{n+s}\sum_{a=0}^{k}\sum_{j=0}^{r}\sum_{l=1}^{s}(-1)^{k+j+a+l-1}\binom{k}{a}\binom{r}{j}\binom{s}{l}2^{l-j-1}\\
&\quad \times\sum_{\substack{i_0 + i_1 + \cdots + i_{s-l} = n+1-l \\ i_0, i_1,\dots, i_{s-l} \ge 1}}\binom{n+s}{i_0} G_{i_1} \cdots G_{i_s}G_{i_0}(j+a \lambda)\bigg\}\,b_{k,\lambda}^{(r)}(x)\\
&=\sum_{k=0}^{n-s}\bigg\{\frac{2^{r}}{\lambda^{k}}\frac{1}{n+s}\sum_{a=k}^{n}\sum_{j=0}^{r}
\sum_{l=1}^{s}\binom{r}{j}\binom{s}{l}\binom{i_{0}}{a}(-2)^{l-j-1}\lambda^{a}S_{2}(a,k)\\
&\quad \times\sum_{\substack{i_0 + i_1 + \cdots + i_{s-l} = n+1-l \\ i_1,\dots, i_{s-l} \ge 1,\,  i_0 \ge a+1}}\binom{n+s}{i_0} G_{i_1} \cdots G_{i_s}G_{i_0-a}(j)\bigg\}\,b_{k,\lambda}^{(r)}(x).
\end{align*}

\section{Conclusion}

In this paper, we were interested in representing any polynomial in terms of the ordered Bell and degenerate ordered Bell polynomials, and more generally of the higher-order ordered Bell and higher-order degenerate ordered Bell polynomials. We were able to derive formulas for such representations with the help of umbral calculus. Further, we illustrated the formulas with some examples. \par
Even though the method adopted in this paper is elementary, they are very useful and powerful. Indeed, as we mentioned in the Section 1, both a variant of Miki's identity and Faber-Pandharipande-Zagier (FPZ) identity follow from the one identity (see \eqref{2A}) that can be derived from a formula (see \eqref{1A}) involving only derivatives and integrals of the given polynomial, while all the other proofs are quite involved. We recall here that the FPZ identity was a conjectural relations between Hodge integrals in Gromov-Witten theory.\par
It is one of our future research projects to continue to find formulas representing polynomials in terms of some specific special polynomials and to apply those in discovering some interesting identities.


\end{document}